\documentclass{SciPost}

\hypersetup{
    colorlinks,
    linkcolor={red!50!black},
    citecolor={blue!50!black},
    urlcolor={blue!80!black}
}

\usepackage[bitstream-charter]{mathdesign}
\usepackage{comment}
\usepackage{amsthm}
\theoremstyle{definition}
\newtheorem{definition}{Definition}[section]
\theoremstyle{plain}
\newtheorem{proposition}[definition]{Proposition}

\newtheorem{corollary}[definition]{Corollary}
\theoremstyle{definition}
\newtheorem{remark}[definition]{Remark}

\DeclareSymbolFont{usualmathcal}{OMS}{cmsy}{m}{n}
\DeclareSymbolFontAlphabet{\mathcal}{usualmathcal}

\newcommand{\RR}{\mathbb{R}}
\newcommand{\TT}{\mathbb{T}}
\newcommand{\Om}{\Omega_\sigma}

\fancypagestyle{SPstyle}{
\fancyhf{}
\lhead{\colorbox{scipostblue}{\bf \color{white} ~SciPost Physics Proceedings }}
\rhead{{\bf \color{scipostdeepblue} ~Submission }}

\fancyfoot[C]{\textbf{\thepage}}
}

\begin{document}

\pagestyle{SPstyle}

\begin{center}{\Large \textbf{\color{scipostdeepblue}{
Symmetry, boundary conditions, and the Noether correspondence of the Zakharov--Kuznetsov equation\\
}}}\end{center}

\begin{center}\textbf{
Yuichiro Shimazaki
}\end{center}

\begin{center}
Department of Physics and Astronomy, Tokyo University of Science, Noda, Chiba 278-8510, Japan
\\[\baselineskip]
\href{mailto:shimazakitus@gmail.com}{\small shimazakitus@gmail.com}
\end{center}

\definecolor{palegray}{gray}{0.95}
\begin{center}
\colorbox{palegray}{
  \begin{tabular}{ll}
 \hspace{-11pt} \begin{minipage}{0.18\textwidth}
    \includegraphics[width=28mm]{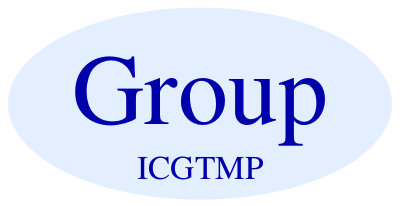}
  \end{minipage}
  &
  \begin{minipage}{0.75\textwidth}
    \begin{center}
    {\it 36th International Colloquium on Group Theoretical Methods in Physics}\\
    {\it Valladolid, 13-17 July 2026} \\
    \doi{10.21468/SciPostPhysProc.?}\\
    \end{center}
  \end{minipage}
\end{tabular}
}
\end{center}

\section*{\color{scipostdeepblue}{Abstract}}
\textbf{\boldmath{%
Conserved quantities of nonlinear dispersive wave equations are typically derived on the whole plane, whereas numerical simulations evaluate them on a torus.
We examine how compactifying one or both spatial directions affects the Noether correspondences of the Zakharov--Kuznetsov equation, obtained from its potential Lagrangian, on the family of domains $\{\RR,\TT\}^2$.
The outcome is a three-valued classification: intrinsic, chart-dependent, or defective.
Only the two first moments, whose characteristics carry an undifferentiated coordinate weight, are affected.
A moment can become defective only if its weighted direction is compact and the mass flux crosses the seam of the chosen fundamental domain; for the Zakharov--Kuznetsov equation this happens exactly for the longitudinal moment when $x$ is compact.
The status of each moment depends only on the boundary condition along its weighted direction.
When $x$ is compact, the Galilean boost still maps solutions to solutions, but no global Noether charge corresponds to it.
}}

\vspace{\baselineskip}

\noindent\textcolor{white!90!black}{
\fbox{\parbox{0.975\linewidth}{
\textcolor{white!40!black}{\begin{tabular}{lr}
  \begin{minipage}{0.6\textwidth}
    {\small Copyright attribution to authors. \newline
    This work is a submission to SciPost Phys. Proc. \newline
    License information to appear upon publication. \newline
    Publication information to appear upon publication.}
  \end{minipage} & \begin{minipage}{0.4\textwidth}
    {\small Received Date \newline Accepted Date \newline Published Date}
  \end{minipage}
\end{tabular}}
}}
}


\vspace{10pt}
\noindent\rule{\textwidth}{1pt}
\tableofcontents
\noindent\rule{\textwidth}{1pt}
\vspace{10pt}

\section{Introduction}
\label{sec:intro}

Conserved quantities of nonlinear dispersive partial differential equations govern the qualitative dynamics of wave packets and serve as diagnostics in numerical simulations.
Analytically, they are usually established on the unbounded plane $\RR^2$, where boundary terms from integration by parts vanish.
Numerical simulations, in contrast, are routinely performed on a torus $\TT^2$, where Fourier pseudospectral methods impose periodic boundary conditions~\cite{Gottlieb1977, Trefethen2000}.
For invariants whose densities are built from the field and its derivatives, the difference is harmless.
For invariants whose densities carry an explicit coordinate weight, it is not: a coordinate is not a function on a circle, so such an invariant must first be defined on a chosen fundamental domain, and its conservation must then be checked directly.

The Zakharov--Kuznetsov (ZK) equation, derived for nonlinear ion-acoustic waves in a magnetized plasma~\cite{Sagdeev1969, Zakharov1974}, is a two-dimensional non-integrable extension of the Korteweg--de Vries (KdV) equation~\cite{Zabusky1965, Miura1968}.
Besides mass, momentum, and energy, it conserves two first moments, $I_{4x}$ and $I_{4y}$~\cite{Zakharov1974}.
They fix the motion of the center of mass and are used to extract propagation speeds and to analyze soliton stability~\cite{Kuznetsov1986, Iwasaki1990, Zakharov2012, Kuznetsov2018}.
Their densities carry the coordinate weights $x$ and $y$, which makes them the natural objects in which a dependence on the domain can appear.
The Cauchy problem of the ZK equation has been studied on the plane and on cylinders~\cite{Faminskii1995, Linares2009, Linares2011, Molinet2015}.

On $\TT^2$, Ref.~\cite{Sawado2026arxiv} integrated the local conservation law directly, obtained a closed-form expression for the boundary defect in $I_{4x}$, and introduced a modified quantity $I_{4x}^{\mathrm{mod}}$ that removes the observed drift.
On $\RR^2$, Adem and Muatjetjeja~\cite{Adem2015} formulated a potential Lagrangian for the ZK equation and classified its Noether point symmetries, and Mothibi and Khalique~\cite{Mothibi2015} analyzed the Lie symmetries of generalized ZK equations.
These classifications yield five symmetry--invariant correspondences, including an affine shift of the potential producing $I_{4y}$ and a Galilean boost producing $I_{4x}$.

The present work follows these correspondences from $\RR^2$ to the domain family $\{\RR,\TT\}^2$, in which each spatial direction is independently non-compact or periodic.
We work within the Noether framework~\cite{Olver1993, Bluman2010}.
The symmetry viewpoint complements the direct computation of Ref.~\cite{Sawado2026arxiv}: it identifies in advance which correspondences can be affected, namely those whose characteristics carry an undifferentiated coordinate weight, and it treats all four domains at once.
We find that the fate of each correspondence is intrinsic, chart-dependent, or defective.
A first moment can become defective only if (1) its weighted direction is compact and (2) the mass flux crosses the seam of the chosen fundamental domain, that is, the mass flux has a component along the weighted direction; a defect then occurs whenever the seam flux does not vanish identically.
If only condition (1) holds, the moment is chart-dependent but strictly conserved.
For the ZK equation, the mass flux points in the $x$-direction, so $I_{4y}$ is never defective, while $I_{4x}$ is defective exactly when $x$ is compact.
The status of each moment depends only on the boundary condition along its own weighted direction, so that $I_{4x}$ is defective on the cylinder $\TT_x\times\RR_y$ but conserved on the cylinder $\RR_x\times\TT_y$, although both have infinite area.

The Galilean boost maps solutions to solutions on all four domains, but when $x$ is compact no global Noether charge corresponds to it.
In this setting the modified invariant of Ref.~\cite{Sawado2026arxiv} is not a recovered charge; we show that it is not a function of the current state and interpret it as a numerical diagnostic.

This paper is organized as follows.
Section~\ref{sec:zk} defines the equation, its invariants, and the local conservation laws.
Section~\ref{sec:symmetry} reviews the potential formulation and the five correspondences on $\RR^2$.
Section~\ref{sec:domains} determines the fate of each correspondence on $\{\RR,\TT\}^2$.
Section~\ref{sec:criterion} formulates the selection rule, and Sec.~\ref{sec:modified} discusses the modified invariant.
Appendix~\ref{app:lagrangian} compiles the Noether identities, and Appendix~\ref{app:defect} the domain-dependent integrations.

\section{The Zakharov--Kuznetsov equation and its invariants}
\label{sec:zk}

\subsection{Governing equation and setup}
\label{sec:zk_equation}

We consider the Zakharov--Kuznetsov (ZK) equation in the form
\begin{equation}
u_t + \alpha\,u\,u_x + \partial_x\left(\nabla^2u\right) = 0,
\qquad \nabla^2=\partial_x^2+\partial_y^2 .
\label{eq:zk}
\end{equation}
The parameter $\alpha$ can be eliminated by rescaling, but we keep it general.
We assume $\alpha\neq0$ where the Galilean boost is discussed.
When periodic boundary conditions are imposed, the periods along $x$ and $y$ are denoted by $P_x$ and $P_y$.

\subsection{Conserved quantities}
\label{sec:zk_invariants}

The ZK equation admits three scalar conserved quantities and two first moments~\cite{Zakharov1974}:
\begin{equation}
I_1 := \iint u\,dx\,dy ,
\qquad
I_2 := \frac{1}{2}\iint u^2\,dx\,dy ,
\qquad
I_3 := \iint\left(\frac{1}{2}|\nabla u|^2-\frac{\alpha}{6}u^3\right)\,dx\,dy,
\label{eq:scalars}
\end{equation}
\begin{equation}
I_{4x} := \iint x\,u\,dx\,dy \;-\; \alpha\,I_2\,t ,
\qquad
I_{4y} := \iint y\,u\,dx\,dy .
\label{eq:moments}
\end{equation}
The factor $\tfrac{1}{2}$ in $I_2$ matches the Noether density of Sec.~\ref{sec:symmetry}, and the subtraction $\alpha I_2t$ coincides with that of Ref.~\cite{Sawado2026arxiv}.
In addition, the equation conserves the family $\iint f(y)\,u$, which contains $I_1$ and $I_{4y}$.

\subsection{Local conservation laws and flux geometry}
\label{sec:zk_flux}

Equation~\eqref{eq:zk} can be written in the local conservation form:
\begin{equation}
\partial_t u + \partial_x F = 0,
\qquad F := \frac{\alpha}{2}u^2+\nabla^2u .
\label{eq:mass_law}
\end{equation}
The mass flux $(F,0)$ has no $y$-component.
Therefore, integrating Eq.~\eqref{eq:mass_law} in $x$ gives $\partial_t\mu=0$ pointwise in $y$, where $\mu(t,y):=\int u\,dx$ is the $y$-line density.
Ref.~\cite{Sawado2026arxiv} noted this conservation as a remnant of the one-dimensional KdV structure.

The absence of a $y$-flux is a property of the mass conservation law, not of every conservation law of the equation.
For instance, the momentum law carries a genuine $y$-flux (Appendix~\ref{app:lagrangian}).
For the first moments, however, only the mass flux matters.
Their densities are coordinate-weighted versions of the mass density, and in the corresponding identities the coordinate weight multiplies only the mass flux (Sec.~\ref{sec:criterion}).

\section{Potential formulation and Noether correspondences on \texorpdfstring{$\RR^2$}{R2}}
\label{sec:symmetry}

\subsection{Potential formulation and Lagrangian density}
\label{sec:symmetry_potential}

Equation~\eqref{eq:zk} is not variational in its original form.
Introducing the potential $v$ via $u=v_x$ yields the potential form of the equation
\begin{equation}
v_{tx}+\alpha v_xv_{xx}+v_{xxxx}+v_{xxyy}=0 ,
\label{eq:potential_eq}
\end{equation}
which is the Euler--Lagrange equation $\mathcal{E}(v)=0$ of the Lagrangian density
\begin{equation}
\mathcal{L}[v]=-\frac{1}{2}v_tv_x-\frac{\alpha}{6}v_x^3+\frac{1}{2}v_{xx}^2+\frac{1}{2}v_{xy}^2 .
\label{eq:lagrangian}
\end{equation}
Since $\mathcal{L}$ depends on $v$ only through its derivatives, transformations of $v$ that leave $u$ invariant can be symmetries of Eq.~\eqref{eq:potential_eq}.
These transformations act trivially on Eq.~\eqref{eq:zk}.
Two of the five correspondences below are of this kind.

\subsection{Variational framework and periodicity of densities}
\label{sec:symmetry_framework}

A characteristic $W$ defines a variational symmetry if $\delta_W\mathcal{L}=D_tB^t+D_xB^x+D_yB^y$ for some $(B^t,B^x,B^y)$, and then Noether's theorem provides the local conservation law
\begin{equation}
W\mathcal{E}(v)=D_tT^t+D_xT^x+D_yT^y .
\label{eq:noether_identity}
\end{equation}
For the Lagrangian~\eqref{eq:lagrangian}, a conserved local functional $G=\iint T^t$ admits the characteristic $W=\delta G/\delta u$, the variational derivative with respect to $u$; all five correspondences below are of this form (Table~\ref{tab:currents}).
The induced transformation of the physical field is $\delta u=\varepsilon\,\partial_xW$.

The potential $v=\int u\,dx$ is periodic in $x$ only if $\mu\equiv0$.
We therefore test periodicity not on $v$ but on the conserved density.
We ask whether $T^t$ defines a functional on the space of periodic fields.
This question is independent of the gauge freedom in $v$, and it reduces to whether the coordinate weights in the characteristics define global functions on the domain.

\subsection{Five Noether correspondences on \texorpdfstring{$\RR^2$}{R2}}
\label{sec:symmetry_five}

The five correspondences of Refs.~\cite{Adem2015,Mothibi2015} are summarized in Table~\ref{tab:currents}, with the characteristics chosen so that $W=\delta G/\delta u$.
Appendix~\ref{app:lagrangian} compiles the associated identities.
Two of them require discussion.

\paragraph{$I_{4y}$ does not arise from $y$-translation.}
The characteristic of $y$-translation is $v_y$, not $y$.
Instead, $I_{4y}$ originates from the affine gauge shift $v\mapsto v+\varepsilon y$, which leaves $u=v_x$ unchanged.
Thus it is a Noether charge associated with a redundancy of the potential representation rather than with a transformation of $u$.
More generally, $v\mapsto v+\varepsilon f(y)$ is a variational symmetry for any $f$; $f\equiv1$ gives $I_1$, and $f(y)=y$ gives $I_{4y}$.
By Corollary~\ref{cor:family}, this whole family remains conserved on every domain.

\paragraph{The Galilean boost and its function-space status.}
The finite boost acting on the potential is
\begin{equation}
v_c(t,x,y)=v(t,x-ct,y)+\frac{c}{\alpha}x-\frac{c^2}{2\alpha}t ,
\label{eq:boost_potential}
\end{equation}
which induces the following transformation of the physical field as
\begin{equation}
(\mathcal{G}_cu)(t,x,y)=u(t,x-ct,y)+\frac{c}{\alpha} .
\label{eq:boost_u}
\end{equation}
The translation and the constant shift are inseparable.
The translation generates $-cu_x$, which is cancelled by the $cu_x$ term induced by the constant shift in the nonlinearity.
The constant $c/\alpha$ reappears in the characteristic $W_{4x}=x/\alpha-tu$ as the undifferentiated weight $x/\alpha$.

As a transformation of the natural phase space, the boost is well defined only when both spatial directions are compact.
On $\RR^2$ the shift $u\mapsto u+c/\alpha$ destroys decay and preserves neither $L^2$ nor $H^1$, so $\mathcal{G}_c$ is only a formal transformation.
The same holds on the cylinders $\TT_x\times\RR_y$ and $\RR_x\times\TT_y$, because a nonzero constant does not decay along the non-compact direction.
Only on $\TT^2$ is $\mathcal{G}_c$ a genuine group action on the phase space.
As shown in Sec.~\ref{sec:domains}, compactifying $x$ removes the corresponding Noether charge; on $\TT^2$ the boost is therefore a genuine symmetry of the phase space without a charge.

\subsection{Characteristics with explicit coordinate weights}
\label{sec:symmetry_weights}

Only two characteristics contain an undifferentiated coordinate weight, namely $W_{4y}=y$ and $W_{4x}=x/\alpha-tu$.
The other three depend only on $u$ and its derivatives; in particular, with the gauge choice of Appendix~\ref{app:lagrangian}, $W_3$ equals $-F$.
A characteristic without coordinate weight is periodic along a compact factor and decays along non-compact directions, so it defines a well-behaved functional on every domain of the family.
The analysis below therefore concerns only $I_{4x}$ and $I_{4y}$.

\begin{table}[htbp]
\centering
\caption{Symmetry--invariant correspondences on $\RR^2$, with characteristics oriented so that $W=\delta G/\delta u$. The last column gives the undifferentiated coordinate weight of the characteristic.}
\label{tab:currents}
\begin{tabular}{llcc}
\hline\hline
$W$ & symmetry & invariant & linear weight \\
\hline
$W_1=1$ & constant shift of $v$ & $I_1$ & --- \\
$W_2=v_x$ & translation in $x$ & $I_2$ & --- \\
$W_3=v_t$ & translation in $t$ & $I_3$ & --- \\
$W_{4y}=y$ & affine shift of $v$ in $y$ & $I_{4y}$ & $y$ \\
$W_{4x}=x/\alpha-tv_x$ & Galilean boost & $I_{4x}/\alpha$ & $x$ \\
\hline\hline
\end{tabular}
\end{table}

\section{Noether correspondences on the domain family \texorpdfstring{$\{\RR,\TT\}^2$}{\{R,T\}2}}
\label{sec:domains}

\subsection{Domain family and three-valued classification}
\label{sec:domains_setup}

Let $\Om=\Omega_x\times\Omega_y$ denote the spatial domain, where $\sigma=(\sigma_x,\sigma_y)\in\{\RR,\TT\}^2$ specifies whether each factor is non-compact or compact, with $\TT_\xi=\RR/P_\xi\mathbb{Z}$.
Along a compact direction, the fields are periodic; along a non-compact direction, they decay rapidly.

On a compact direction $\xi$, a coordinate weight is defined only after a fundamental domain $[s,s+P_\xi)$ has been chosen.
We refer to the choice of $s$ as a chart.
Changing the origin, $\xi\mapsto\xi-a$, changes a first moment by $-aI_1$, as on $\RR$.
Since $I_1$ is conserved, the origin is immaterial, and the essential choice is the seam.

\begin{definition}[Three-valued classification]
\label{def:status}
Let $I$ be a local functional of $(t,u(t))$, defined in each chart when coordinate weights on compact directions are involved.
It is
\begin{enumerate}
\item \emph{intrinsic} if it does not depend on the chart and is conserved along every solution;
\item \emph{chart-dependent} if it is conserved along every solution in every chart, but its value depends on the chart;
\item \emph{defective} if, in some chart, it fails to be conserved along some solution.
\end{enumerate}
\end{definition}

The restriction to local functionals of the instantaneous state is essential.
For any functional $G$ and solution map $\Phi_t$, the functional $u\mapsto G(\Phi_{-t}u)$ is trivially conserved.
Allowing such functionals would make every quantity conserved.
Section~\ref{sec:modified} shows that the modified invariant of Ref.~\cite{Sawado2026arxiv} falls outside Definition~\ref{def:status} for this reason.

\subsection{Intrinsic status of scalar invariants \texorpdfstring{$I_1$, $I_2$, and $I_3$}{I1, I2, and I3}}
\label{sec:domains_scalars}

The characteristics $W_1=1$, $W_2=u$, and $W_3=-F$ contain no coordinate weight.
The Noether identity~\eqref{eq:noether_identity} integrates to a conservation law on every $\Om$, the boundary terms vanishing by periodicity or decay.
Hence $I_1$, $I_2$, and $I_3$ are intrinsic on the whole family, and the domain dependence studied below is a property of the first moments alone.

\subsection{Exact conservation and chart dependence of \texorpdfstring{$I_{4y}$}{I4y}}
\label{sec:domains_I4y}

When $y$ is compact, the weight $W_{4y}=y$ is not a global function on $\Om$.
Nevertheless, $I_{4y}$ is exactly conserved on all four domains, as the identity
\begin{equation}
y\,\mathcal{E}(v)=D_t(yu)+D_x(yF)
\label{eq:cl_I4y}
\end{equation}
shows.
It contains no $D_y$ term: the weight $y$ multiplies the mass flux inside an $x$-derivative and is never evaluated on a $y$-boundary.
Integration in $x$ alone closes the balance by periodicity or decay (Proposition~\ref{prop:linedensity}).

When $y$ is compact, the value of $I_{4y}$ depends on the chart through a conserved but generally nonzero amount (Proposition~\ref{prop:chartshift}).
Its status is therefore chart-dependent, while its conservation is untouched.
When $y$ is non-compact, $y$ is a global coordinate and $I_{4y}$ is intrinsic.
Chart dependence is a property of the linear representative: for compact $y$, the functional $\iint e^{2\pi iy/P_y}u$, a member of the family in Corollary~\ref{cor:family}, is intrinsic and encodes the transverse position of the center of mass on the circle.

\subsection{Breakdown of the Noether correspondence for \texorpdfstring{$I_{4x}$}{I4x} on compact \texorpdfstring{$x$}{x}}
\label{sec:domains_I4x}

For $I_{4x}$, both ingredients of a defect are present when $x$ is compact: the weight $x/\alpha$ is not a global function, and $x$ is the direction of the mass flux.
Three statements characterize the situation.

\paragraph{Non-existence of a global charge when $x$ is compact.}
Suppose that a conserved local functional $G$ of periodic fields corresponded to the boost, in the sense that its characteristic $\delta G/\delta u$ generates the same transformation of $u$ as $W_{4x}$.
Since $G$ is a functional of periodic fields, $\delta G/\delta u$ is $P_x$-periodic in $x$.
Because $\delta u=\varepsilon\,\partial_xW$, two characteristics generate the same transformation of $u$ only if they differ by a function $h(t,y)$, which is precisely the gauge freedom of $v$.
Hence $x/\alpha-tu+h(t,y)$ would have to be periodic in $x$ for some $h$, which is impossible, since its $x$-derivative $1/\alpha-tu_x$ has integral $P_x/\alpha\neq0$ over one period.
The same obstruction is also visible in the line density.
The boost~\eqref{eq:boost_u} raises the conserved line density $\mu$ by $cP_x/\alpha$, whereas the transformation $\delta u=\varepsilon\,\partial_xW$ generated by any $x$-periodic characteristic leaves $\mu$ unchanged.
The argument uses only the compactness of $x$ and is independent of the boundary condition in $y$.

\paragraph{Existence when $x$ is non-compact.}
When $x$ is non-compact, $x$ is a global coordinate, and Proposition~\ref{prop:I4xR} shows directly that $I_{4x}$ is finite, conserved, and chart-independent on $\RR_x\times\RR_y$ and on $\RR_x\times\TT_y$.

\paragraph{Boundary defect when $x$ is compact.}
In a fixed chart, the drift rate of $I_{4x}$ is given in closed form by Proposition~\ref{prop:defect}: it equals $-P_x$ times the net mass flux through the seam.
This flux is not identically zero in general, so $I_{4x}$ is defective.

\subsection{Persistence of symmetry versus breakdown of the Noether correspondence}
\label{sec:domains_severed}

The boost~\eqref{eq:boost_u} maps periodic solutions to periodic solutions, with the translation understood modulo $P_x$, and more generally maps solutions to solutions on every $\Om$.
On $\TT^2$, it is moreover a genuine group action on the phase space (Sec.~\ref{sec:symmetry_five}).
What fails when $x$ is compact is the lift to a Noether charge: the generator $X=t\,\partial_x+(x/\alpha)\,\partial_v$ contains the non-periodic coefficient $x/\alpha$, and no choice of gauge makes its characteristic periodic.
The correspondence between the symmetry and a conserved charge is severed, while the symmetry of the equation itself persists.

\subsection{Classification summary across domains}
\label{sec:domains_table}

\begin{table}[htbp]
\centering
\caption{Classification of each invariant across the domain family $\{\RR,\TT\}^2$.}
\label{tab:fate}
\begin{tabular}{lccc}
\hline\hline
$(\sigma_x,\sigma_y)$ & $I_1,I_2,I_3$ & $I_{4y}$ & $I_{4x}$ \\
\hline
$(\RR,\RR)$ & intrinsic & intrinsic & intrinsic \\
$(\RR,\TT)$ & intrinsic & chart-dependent & intrinsic \\
$(\TT,\RR)$ & intrinsic & intrinsic & defective \\
$(\TT,\TT)$ & intrinsic & chart-dependent & defective \\
\hline\hline
\end{tabular}
\end{table}

On $(\RR,\TT)$ one direction is compact, yet $I_{4x}$ remains intrinsic.
On $(\TT,\RR)$ the domain has infinite area, yet $I_{4x}$ is defective.
The relevant factor is thus not the size of the domain.

The fully periodic case $\TT^2$ agrees with the numerical results of Ref.~\cite{Sawado2026arxiv}, where the measured drift of $I_{4x}$ is accounted for by the boundary flux.
The statements on the mixed domains rest on the analytical results of Appendix~\ref{app:defect}.

\section{Selection rule}
\label{sec:criterion}

\subsection{Formulation of the selection rule}
\label{sec:criterion_statement}

For a first moment with weight $\xi$, the relevant identity is obtained by multiplying the mass conservation law~\eqref{eq:mass_law} by $\xi$:
\begin{equation}
\xi\,\mathcal{E}(v)=D_t(\xi u)+D_x(\xi F)-\delta_{\xi x}F ,
\label{eq:weighted_mass}
\end{equation}
which reproduces Eq.~\eqref{eq:cl_I4y} for $\xi=y$ and the weighted part of Eq.~\eqref{eq:cl_I4x} for $\xi=x$.
In a chart with seam $s$, the weight $\xi$ jumps by $-P_\xi$ across the seam.
This jump contributes to the balance only when the mass flux crosses the seam, that is, when the weight appears inside a derivative along $\xi$.
This leads to the following rule.

\begin{proposition}[Selection rule]
\label{prop:rule}
Consider a first moment with a coordinate weight $\xi$, where $\xi$ is $x$ or $y$.
\begin{enumerate}
\item It is intrinsic if $\xi$ is non-compact.
\item If $\xi$ is compact, it can be defective only if there is a mass flux in the $\xi$ direction; in that case its drift rate in a chart with seam $s$ is $-P_\xi$ times the net mass flux through the seam.
\item If $\xi$ is compact and the mass flux has no component along $\xi$, it is chart-dependent.
\end{enumerate}
It is defective if and only if $\xi$ is compact, the mass flux has a component along $\xi$, and the net mass flux through the seam does not vanish identically.
\end{proposition}

For the ZK equation, the mass flux is $(F,0)$ and the seam flux is not identically zero (Proposition~\ref{prop:defect}), which yields Table~\ref{tab:fate}.
The non-degeneracy condition is not automatic.

\begin{remark}
\label{rem:nondeg}
The linear equation $u_t+u_{xyy}=0$ is the Euler--Lagrange equation of the Lagrangian $\mathcal{L}=-\frac12v_tv_x+\frac12v_{xy}^2$.
The mass flux $(u_{yy},0)$ crosses a seam in $x$, but the net seam flux is the integral of $u_{yy}$ over $\Omega_y$ and vanishes identically.
The $x$-moment $\iint x\,u$ is therefore chart-dependent, not defective, on domains with compact $x$.
\end{remark}

Table~\ref{tab:fate} also exhibits a factorization.
The structural status of each moment depends only on the boundary condition along its own weighted direction.
This holds for linear weights in general.
The lift of $\xi$ jumps only across seams transverse to $\xi$ and is continuous across the seams of the other direction, where the flux terms of the identity are periodic or decaying and therefore cancel.
Consequently, the boundary condition in the transverse direction enters only through the domain of the transverse integration in the seam flux.

\subsection{Role of the four-domain family}
\label{sec:criterion_counterexample}

On $\TT^2$ alone, the two moments already behave differently: both weighted directions are compact, but only $I_{4x}$ is defective.
This shows that compactness of the weighted direction does not suffice for a defect, and that the direction of the mass flux is a second, independent factor.
What the two mixed domains add is the factorization of Proposition~\ref{prop:rule}.
On $\TT_x\times\RR_y$ and $\RR_x\times\TT_y$, each moment behaves as on the domain on which only its own direction is compactified, independently of the transverse boundary condition.
This cannot be seen from the comparison of $\RR^2$ and $\TT^2$ alone, in which both directions change at once.

\subsection{Structural obstruction and dynamical mechanism}
\label{sec:criterion_layers}

The absence of a global charge does not by itself produce an observed drift.
The right-hand side of Eq.~\eqref{eq:defect_general} vanishes whenever the net mass flux through the seam vanishes.
Observing a drift requires a dynamical mechanism, typically dispersive radiation that reaches the seam and carries mass across it.
The structural statement (non-existence of the charge) and the dynamical one (emergence of the drift) should be kept apart.

The two statements also behave differently under enlargement of the periodic box.
The non-existence argument of Sec.~\ref{sec:domains_I4x} holds for every period $P_x$.
For localized initial data, however, the seam can be placed far from the bulk of the solution, so that the seam flux, and with it the drift, remains small up to a given time.
Finite-time observables approach their values on $\RR^2$ as the box grows, although the charge is not restored for any finite period.

\section{Interpretation of the modified invariant}
\label{sec:modified}

For compact $x$, Ref.~\cite{Sawado2026arxiv} introduced the modified quantity
\begin{equation}
I_{4x}^{\mathrm{mod}}(t)
:= \widetilde I_{4x}^{\,(x_0)}(t)-\int_{0}^{t}\mathcal{A}(\tau)\,d\tau,
\qquad
\mathcal{A}(t):=-P_x\!\int_{\Omega_y}F(t,x_0,y)\,dy ,
\label{eq:modified}
\end{equation}
and showed that it removes the observed drift to numerical precision.
Let
\begin{equation}
N_{x_0}(t):=\int_0^t\!\!\int_{\Omega_y}F(\tau,x_0,y)\,dy\,d\tau
\label{eq:netmass}
\end{equation}
be the net mass transported across the seam $x=x_0$ in the positive $x$-direction up to time $t$.

\begin{proposition}
\label{prop:modified}
Let $\Omega_x=\TT$.
Then
\begin{equation}
I_{4x}^{\mathrm{mod}}(t)=\widetilde I_{4x}^{\,(x_0)}(t)+P_xN_{x_0}(t)=\widetilde I_{4x}^{\,(x_0)}(0)
\qquad\text{for all }t .
\label{eq:modified_is_initial}
\end{equation}
Moreover, there is no functional $\Phi$ on the space of periodic fields such that $\Phi[u(t)]=I_{4x}^{\mathrm{mod}}(t)$ for every solution and every $t$.
\end{proposition}

\begin{proof}
The first statement follows from Eqs.~\eqref{eq:modified} and~\eqref{eq:netmass} and Proposition~\ref{prop:defect}.
For the second, at $t=0$, such a $\Phi$ would have to coincide with $\widetilde I_{4x}^{\,(x_0)}|{t=0}=\iint x,u$ on all periodic fields.
But $\iint x\,u$ changes along solutions at the rate $\alpha I_2-P_x\int_{\Omega_y}F(t,x_0,y)\,dy$, which is nonzero for suitable data, for example for small-amplitude data supported away from the seam when $\alpha\neq0$.
\end{proof}

Equation~\eqref{eq:modified_is_initial} has a simple meaning.
Mass that crosses the seam in the positive direction re-enters the fundamental domain with its coordinate reduced by $P_x$, and the term $P_x\,N_{x_0}$ restores this contribution.
The modified quantity is thus constant, but its constant value is simply the initial value of the chart integral.
It is not a recovered Noether charge, which does not exist when $x$ is compact (Sec.~\ref{sec:domains_I4x}).
And by Proposition~\ref{prop:modified} it is not a function of the current state, so it lies outside Definition~\ref{def:status}.
It remains chart-dependent through $x_0$.
Its value lies in its use as a diagnostic.
The constancy verifies that a simulation resolves the balance of Proposition~\ref{prop:defect}.

The four-domain framework also fixes the scope of this construction.
On $(\RR,\RR)$ and $(\RR,\TT)$, Proposition~\ref{prop:I4xR} gives $\mathcal{A}\equiv0$ and the correction is unnecessary.
It is required on $(\TT,\RR)$ and $(\TT,\TT)$, where Eq.~\eqref{eq:modified} applies in the same form, the boundary condition along $y$ entering only through the integration domain of $\mathcal{A}$.

\section{Conclusion}
\label{sec:conclusion}

We followed the five Noether correspondences of the Zakharov--Kuznetsov equation across the domain family ${\RR,\TT}^2$.
The three scalar correspondences are insensitive to the choice of domain.
The transverse moment $I_{4y}$ is conserved on all four domains and becomes chart-dependent when $y$ is compact.
The longitudinal moment $I_{4x}$ has no global Noether charge when $x$ is compact and is then defective, at a rate given by the net mass flux through the seam.

The outcome is governed by the compactness of the weighted direction together with the direction of the mass flux, supplemented by a non-degeneracy condition on the seam flux.
The status of each moment depends only on the boundary condition along its own weighted direction.
The two mixed domains exhibit this factorization, which a comparison of $\RR^2$ and $\TT^2$ alone cannot.

The Galilean boost maps solutions to solutions on every domain and is a genuine group action on the phase space of $\TT^2$, yet compactifying $x$ removes its Noether charge.
The modified invariant of Ref.~\cite{Sawado2026arxiv} is not a recovered charge and not a function of the current state.
It records the net mass transported across the seam and serves as a numerical diagnostic.

\section*{Acknowledgements}
The author is grateful to Nobuyuki Sawado, Atsushi Nakamula, Kouichi Toda, Naruhiko Aizawa and E. A. Kuznetsov for meaningful discussions and helpful comments.
The author also thanks the organizers of the 36th International Colloquium on Group Theoretical Methods in Physics (Group 36).

\paragraph{Funding information}
This work was supported in part by JSPS KAKENHI Grant Number JP23K02794.

\begin{appendix}
\numberwithin{equation}{section}

\section{Lagrangian structure and Noether identities}
\label{app:lagrangian}

This appendix records the Euler--Lagrange equation and the Noether identities in the jet space; they hold pointwise on every domain $\Om$.

\subsection{Euler--Lagrange equation}

For the Lagrangian~\eqref{eq:lagrangian}, the Euler--Lagrange operator is
\begin{equation}
\mathcal{E}(v)
= -D_t\frac{\partial\mathcal{L}}{\partial v_t}
-D_x\frac{\partial\mathcal{L}}{\partial v_x}
+D_x^2\frac{\partial\mathcal{L}}{\partial v_{xx}}
+D_xD_y\frac{\partial\mathcal{L}}{\partial v_{xy}} .
\label{eq:el_operator}
\end{equation}
With $\partial\mathcal{L}/\partial v_t=-\tfrac12v_x$, $\partial\mathcal{L}/\partial v_x=-\tfrac12v_t-\tfrac{\alpha}{2}v_x^2$, $\partial\mathcal{L}/\partial v_{xx}=v_{xx}$, and $\partial\mathcal{L}/\partial v_{xy}=v_{xy}$, this gives
\begin{equation*}
\mathcal{E}(v) = v_{tx}+\alpha v_xv_{xx}+v_{xxxx}+v_{xxyy}=u_t+\partial_xF ,
\end{equation*}
which reproduces Eq.~\eqref{eq:potential_eq}.
Since $\mathcal{L}$ depends on $v$ only through $v_t$, $v_x$, $v_{xx}$, and $v_{xy}$, the characteristics $W=f(y)$ satisfy $\delta_W\mathcal{L}=0$ identically.

\subsection{Gauge fixing}

Because Eq.~\eqref{eq:potential_eq} is a total $x$-derivative, the quantity $v_t+\tfrac{\alpha}{2}v_x^2+v_{xxx}+v_{xyy}$ is independent of $x$.
A gauge shift $v \mapsto v + g(t,y)$ leaves $u=v_x$ unchanged and sets
\begin{equation}
v_t=-\left(\frac{\alpha}{2}u^2+\nabla^2u\right)=-F .
\label{eq:vt_is_F}
\end{equation}
Hence $W_3=v_t=-F$ is periodic whenever $u$ is periodic.

\subsection{Standard Noether identities}
\label{app:lagrangian_routine}

Contracting $\mathcal{E}(v)$ with the scalar characteristics yields
\begin{align}
W_1=1: &\quad \mathcal{E}(v) = D_tu + D_xF,
\label{eq:cl_mass} \\
W_2=u: &\quad u\,\mathcal{E}(v) = D_t\left(\frac12u^2\right) + D_x\Phi_x + D_y\Phi_y,
\label{eq:cl_mom} \\
W_3=-F: &\quad -F\,\mathcal{E}(v) = D_t h + D_x \Psi_x + D_y \Psi_y,
\label{eq:cl_ene} \\
W_{4y}=y: &\quad y\,\mathcal{E}(v) = D_t(yu) + D_x(yF) ,
\label{eq:cl_I4y_app}
\end{align}
where $h=\tfrac12|\nabla u|^2-\tfrac{\alpha}{6}u^3$ is the density of $I_3$ and
\begin{equation*}
\Phi_x=\frac{\alpha}{3}u^3+uu_{xx}-\frac{1}{2}u_x^2+uu_{yy}+\frac{1}{2}u_y^2,
\quad
\Phi_y=-u_xu_y,
\quad
\Psi_x=-\frac12F^2-u_xu_t,
\quad
\Psi_y=-u_yu_t .
\end{equation*}
Equation~\eqref{eq:cl_mom} is the momentum law, whose flux has a nonzero $y$-component, whereas Eqs.~\eqref{eq:cl_mass} and~\eqref{eq:cl_I4y_app} contain no $D_y$ term.

\subsection{Galilean boost}

For $W_{4x}=x/\alpha-tu$, decomposing the contraction via the mass and momentum identities yields
\begin{align}
\frac{x}{\alpha}\,\mathcal{E}(v)
&= D_t\left(\frac{x}{\alpha}u\right)+D_x\left(\frac{x}{\alpha}F\right)-\frac{1}{\alpha}F,
\label{eq:boost_part1}\\
-tu\,\mathcal{E}(v)
&= D_t\left(-\frac{t}{2}u^2\right)-D_x(t\Phi_x)-D_y(t\Phi_y)+\frac{1}{2}u^2 .
\label{eq:boost_part2}
\end{align}
The remaining undifferentiated terms combine as
\begin{equation}
-\frac{1}{\alpha}F+\frac{1}{2}u^2
= -\frac{1}{\alpha}\nabla^2u
= -D_x\left(\frac{u_x}{\alpha}\right)-D_y\left(\frac{u_y}{\alpha}\right) .
\label{eq:boost_remainder}
\end{equation}
Summing Eqs.~\eqref{eq:boost_part1}--\eqref{eq:boost_remainder} gives
\begin{equation}
W_{4x}\,\mathcal{E}(v)
= D_t\left(\frac{x}{\alpha}u-\frac{t}{2}u^2\right)
+ D_x\left(\frac{x}{\alpha}F-t\Phi_x-\frac{u_x}{\alpha}\right)
+ D_y\left(-t\Phi_y-\frac{u_y}{\alpha}\right) .
\label{eq:cl_I4x}
\end{equation}
The coordinate weight appears only in the $x$-flux, through the mass flux $(x/\alpha)F$; the $y$-flux carries no coordinate weight.

\section{Moment invariants on the domain family \texorpdfstring{$\{\RR,\TT\}^2$}{\{R,T\}2}}
\label{app:defect}

\subsection{Decay hypotheses and integration identities}
\label{app:defect_hypotheses}

We consider solutions that are smooth, periodic along compact directions, and rapidly decaying along non-compact ones: all derivatives decay faster than any power of the non-compact coordinates, uniformly in the remaining variables and locally uniformly in $t$, together with their time derivatives.
This persistence of decay is assumed rather than proved.
It ensures that all integrals below converge absolutely, that differentiation under the integral sign is allowed, and that
\begin{equation}
\int_{\Omega_x}\partial_xF\,dx=0,
\qquad
\iint_{\Om}\nabla^2u\,dx\,dy=0
\label{eq:flux_lap_int}
\end{equation}
on every domain $\Om$, $\sigma\in\{\RR,\TT\}^2$.

\subsection{Conservation of \texorpdfstring{$y$}{y}-line density}

\begin{proposition}[Ref.~\cite{Sawado2026arxiv}]
\label{prop:linedensity}
The $y$-line density $\mu(t,y):=\int_{\Omega_x}u(t,x,y)\,dx$ satisfies $\partial_t\mu=0$ pointwise in $y$ on every domain $\Om$.
\end{proposition}

\begin{proof}
Integrate Eq.~\eqref{eq:mass_law} in $x$ at fixed $y$ and use Eq.~\eqref{eq:flux_lap_int}.
\end{proof}

\begin{corollary}
\label{cor:family}
For every measurable $f$ on $\Omega_y$ that is bounded (for compact $y$) or polynomially bounded (for non-compact $y$), the functional $\iint_{\Om} f(y)\,u\,dx\,dy$ is conserved on every domain $\Om$.
When $y$ is compact and $f$ is not periodic, the statement holds in every chart.
\end{corollary}

\subsection{Chart dependence of \texorpdfstring{$I_{4y}$}{I4y}}
\label{app:defect_chart_I4y}

When $y$ is compact, a chart with seam $s$ uses the fundamental domain $[s,s+P_y)$, and the moment in this chart is
\begin{equation}
\widetilde I_{4y}^{\,(s)}
:= \int_{\Omega_x}\!\int_{s}^{\,s+P_y}y\,u\,dy\,dx .
\label{eq:I4y_chart}
\end{equation}

\begin{proposition}
\label{prop:chartshift}
On any domain with compact $y$, and for $s<s'<s+P_y$,
\begin{equation*}
\widetilde I_{4y}^{\,(s')}-\widetilde I_{4y}^{\,(s)}
=P_y\int_s^{s'}\mu(y)\,dy .
\end{equation*}
The right-hand side is conserved but does not vanish identically.
\end{proposition}

\begin{proof}
By periodicity, $\frac{d}{ds'}\widetilde I_{4y}^{\,(s')}=(s'+P_y)\mu(s')-s'\mu(s')=P_y\mu(s')$.
Integrating from $s$ to $s'$ gives the result, and conservation follows from Proposition~\ref{prop:linedensity}.
\end{proof}

\subsection{Conservation of \texorpdfstring{$I_{4x}$}{I4x} for non-compact \texorpdfstring{$x$}{x}}
\label{app:defect_I4xR}

\begin{proposition}
\label{prop:I4xR}
Let $\Omega_x=\RR$, with $\Omega_y\in\{\RR,\TT\}$.
Then $I_{4x}$ defined in Eq.~\eqref{eq:moments} is finite, conserved, and chart-independent.
\end{proposition}

\begin{proof}
Finiteness follows from the decay hypotheses.
Multiplying Eq.~\eqref{eq:mass_law} by $x$ and integrating by parts in $x$ yields
\begin{equation*}
\frac{d}{dt}\iint x\,u
= -\int_{\Omega_y}\left[xF\right]_{x=-\infty}^{+\infty}dy+\iint F .
\end{equation*}
The boundary term vanishes by rapid decay.
By Eq.~\eqref{eq:flux_lap_int}, the area integral of $F$ equals $\tfrac{\alpha}{2}\iint u^2=\alpha I_2$.
Since $I_2$ is conserved, $dI_{4x}/dt=0$.
Chart independence is immediate since $x$ is a global coordinate on $\RR$.
\end{proof}

\subsection{Boundary defect rate for compact \texorpdfstring{$x$}{x}}
\label{app:defect_I4xT}

\begin{proposition}
\label{prop:defect}
Let $\Omega_x=\TT$ with fundamental domain $[x_0,x_0+P_x)$, and $\Omega_y\in\{\RR,\TT\}$.
Define $\widetilde I_{4x}^{\,(x_0)}:=\int_{x_0}^{x_0+P_x}\int_{\Omega_y}\left(x\,u-\tfrac{\alpha}{2}t\,u^2\right)dy\,dx$.
Then
\begin{equation}
\frac{d}{dt}\widetilde I_{4x}^{\,(x_0)}
= -P_x\int_{\Omega_y}F(t,x_0,y)\,dy
= -P_x\int_{\Omega_y}\left(\frac{\alpha}{2}u^2+u_{xx}\right)(t,x_0,y)\,dy .
\label{eq:defect_general}
\end{equation}
The right-hand side is $-P_x$ times the net mass flux through the seam, and it does not vanish identically.
\end{proposition}

\begin{proof}
Integrating by parts in $x$ over the fundamental domain gives
\begin{equation*}
\frac{d}{dt}\iint x\,u
= -\int_{\Omega_y}\left[xF\right]_{x=x_0}^{x_0+P_x}dy+\iint F .
\end{equation*}
Since $F$ is $P_x$-periodic in $x$, the boundary term equals $P_x\int_{\Omega_y}F(t,x_0,y)\,dy$.
By Eq.~\eqref{eq:flux_lap_int}, $\iint_{\Om}F\,dx\,dy=\alpha I_2$, which cancels the time derivative of $-\tfrac{\alpha}{2}t\iint_{\Om}u^2$.
The term $u_{yy}$ in $F$ integrates to zero over $\Omega_y$.
For the datum $u_0=\big(1-\cos(2\pi(x-x_0)/P_x)\big)\psi(y)$ with $\int\psi\,dy\neq0$, one has $u_0(x_0,\cdot)=0$ and a seam flux $(2\pi/P_x)^2\int\psi\,dy\neq0$ at $t=0$.
\end{proof}

\end{appendix}

\bibliography{g36_proc_refs}
\bibliographystyle{SciPost_bibstyle}

\end{document}